\documentclass[11pt]{amsart}
\usepackage{mathrsfs}
\usepackage{amssymb}

\usepackage{amscd}
\usepackage{amsmath, amssymb}
\usepackage{amsfonts}
\usepackage[colorlinks,linkcolor=blue,citecolor=blue, pdfstartview=FitH]{hyperref}
\usepackage[english]{babel}
\usepackage[all]{xy}

\newtheorem{theorem}{Theorem}[section]
\newtheorem{lemma}[theorem]{Lemma}

\newtheorem{corollary}[theorem]{Corollary}

\newtheorem{definition}[theorem]{Definition}

\newtheorem{conjecture}[theorem]{Conjecture}
\newtheorem{remark}[theorem]{Remark}

\newcommand{\pa}{\partial}

\numberwithin{equation}{section}

\usepackage{fancyhdr}
\begin{document}
 \setcounter{section}{-1}
\title{Remarks on BCOV invariants and degenerations of Calabi-Yau manifolds}

\makeatletter
\let\uppercasenonmath\@gobble% disables title uppercase
\let\MakeUppercase\relax% disables author uppercase
\let\scshape\relax% disables section smallcaps
\makeatother

\author{Kefeng Liu}
\address{Kefeng Liu, Department of Mathematics,
Capital Normal University, Beijing, 100048, China}
\address{Department of Mathematics, University of California at Los Angeles, California 90095}

\email{liu@math.ucla.edu}

\author{Wei Xia}

\address{Wei Xia, Center of Mathematical Sciences, Zhejiang
University, Hangzhou, 310027, CHINA. } \email{xiaweiwei3@126.com}

\begin{abstract} For a one parameter family of Calabi-Yau threefolds, Green, Griffiths and Kerr have expressed the total singularities in terms of the degrees of Hodge bundles and Euler number of the general fiber. In this paper, we show that the total singularities can be expressed by the sum of asymptotic values of BCOV invariants, studied by Fang, Lu and Yoshikawa. On the other hand, by using Yau's Schwarz lemma, we prove Arakelov type inequalities and Euler number bound for Calabi-Yau family over a compact Riemann surface. \end{abstract}
\maketitle

\begin{section}{Introduction}
    Reidemeister torsion (R-torsion) is an invariant that can distinguish between closed manifolds which are homotopy equivalent but not homeomorphic. Analytic torsion (or Ray-Singer torsion) is an invariant of Riemannian manifolds defined by Ray and Singer~\cite{RS71,RS73} as an analytic analogue of Reidemeister torsion. These two torsions naturally coincide, which is known as the Cheeger-M\"{u}ller theorem~\cite{Mul78,Che79}.\par

     %i.e. it's invariant under topological deformations but not homotopic ones.%

     In~\cite{BCOV93a,BCOV93b}, Bershadsky-Cecotti-Ooguri-Vafa conjectured an equivalence between the physical quantity $F_1$ of a Calabi-Yau threefold and a linear combination of the holomorphic analytic torsions, which now called the BCOV torsions. Motivated by their conjecture,
     Fang, Lu and Yoshikawa~\cite{FLY08} considered a modification of the BCOV torsions, called the BCOV invariants. They have conducted a detailed study of the asymptotic behavior of the BCOV invariants for the mirror quintic threefold. Recently, Yoshikawa~\cite{Yo14,Yo14a} has continued the study of the asymptotic behavior of the BCOV invariants for general degenerations of Calabi-Yau threefolds.\par

     By using the curvature formula for Quillen metrics~\cite{BGS88a,BGS88b,BGS88c}, Bershadsky-Cecotti-Ooguri-Vafa obtained a variational formula for the BCOV torsion of Ricci-flat Calabi-Yau manifolds. Fang-Lu~\cite{FL05} expressed the variation of the BCOV torsion of Ricci-flat Calabi-Yau manifolds as a linear combination of the Weil-Petersson metric and the generalized Hodge metrics. Using this, they managed to show that for primitive (i.e. $\omega_{H^k}=0, \forall k<n$, see Section \ref{Generalized Hodge metrics and Hodge bundles}) Calabi-Yau manifolds of dimension $n$, if the Euler characteristic number satisfy $(-1)^n\chi  < 24$, then there exists no complete curve in $\mathcal{M}$, the moduli space of polarized Calabi-Yau manifolds.\par

     In their paper~\cite{GGK09}, Green, Griffiths and Kerr have studied global enumerative invariants of a polarized variation of Hodge structures
     for Calabi-Yau manifolds. By applying the \textit{Grothendieck-Riemann-Roch} theorem, they can express the total singularities in terms of the degrees of Hodge bundles and Euler number of the general fiber:
     \begin{equation}\label{GGK main 0}
 (\chi+36)\deg\mathcal{H}_e^{3,0}+12 \deg\mathcal{H}_e^{2,1}=2~\sum_{j}\Delta_j+\sum_i\{\chi_2^i-3I_2^i\},
\end{equation}
where $\mathcal{H}_e^{3,0}, \mathcal{H}_e^{2,1}$ are the canonical Deligne's extensions of the Hodge bundles $\mathcal{H}^{3,0}, \mathcal{H}^{2,1}$, for each $j$, $\Delta_j$ is the number of ODP's on the singular fiber, and for each $i$, $\chi_2^i-3I_2^i$ is the topological and intersecting data of the singular fiber. See \eqref{GGK main} for precise definition. As a corollary, they showed that a non-degenerate family of Calabi-Yau threefolds must satisfy $\chi \leq -24$\footnote{In their paper, this inequality is strict which we believe need further explanation. Through email correspondence, Kerr has told us that he agrees with this point.}.\par

     Let $\mathcal{X}$ be a smooth projective variety of dimension $n+1$ and let $S$ be a compact Riemann surface, $f \colon{\mathcal{X}}\to S$ be a surjective, flat holomorphic map with generic fiber Calabi-Yau $n$-fold, and $f^0 \colon{\mathcal{X}}^0 \to S^0$ be the smooth part of $f$. Unless stated otherwise we shall assume all the local monodromy transformations around points $s_i\in S\setminus S^0$ are unipotent. This is our basic assumption throughout this paper. By integrating the current equation due to Fang, Lu and Yoshikawa,  we find that (for $n=3$) the total singularities can be expressed by the sum of asymptotic values of BCOV invariants:
     \begin{theorem}\label{main theorem 1 0}
 Let $E:=S\setminus S^0=\{s_i\}\subset S$ be the discriminant locus of $f$ and $\tau_{\rm BCOV}(t)$ be the BCOV invariant of the fiber $f^{-1}(t)$. Set
 $$a_i
:=
\lim_{t\to s_i}
\frac{\log\tau_{\rm BCOV}(t)}
{\log|t|^{2}}.$$
Then
\begin{equation}\label{0.1}
(\chi+36)\deg\mathcal{H}_e^{3,0}+12 \deg\mathcal{H}_e^{2,1}=12\sum_{i}~a_i.
\end{equation}
\end{theorem}
The interesting feature of this identity is that the left hand side involves only the zero spectrum, or the harmonic space, of the Laplacian operator, while the right hand side is derived from the whole spectrum. As a direct corollary, we see that

\begin{corollary}$$12\sum_{i}a_i\in \mathbb{Z} .$$
\end{corollary}
 Each individual $a_i$ was previously only known to be a rational number~\cite[Th.\,0.1]{Yo14}. In fact, Yoshikawa has obtained an explicit expression for $a_i$, but it is too complicated. Nevertheless, if the singular fibers of $f$ has at most ODP singularities, then $a_i$ has a simple form. Combining this with Theorem \ref{main theorem 1 0}, we get

 \begin{corollary} If the singular fibers of $f$ has at most ODP singularities, then
 \begin{equation}\label{ODP 0}
 (\chi+36)\deg\mathcal{H}_e^{3,0}+12 \deg\mathcal{H}_e^{2,1}=2\sum_{j}\Delta_j,
 \end{equation}
 where $\Delta_j= \#{\rm Sing}\,X_{j}$ is the number of ODP's on $X_j$.
 \end{corollary}

 Note that \eqref{ODP 0} is a special case of \eqref{GGK main 0}. Conversely, by comparing
\eqref{GGK main 0} and \eqref{0.1}, we have the following geometric meaning for $12\sum_{i}a_i$:

\begin{corollary}\label{semi-stable cor 0}
If $f$ is relatively minimal and has at most semi-stable singular fibers, then
\begin{equation}\label{semi-stable formula 0}
  12\sum_{i}a_i=\sum_i\{\chi_2^i-3I_2^i\}.
\end{equation}
\end{corollary}

By comparing the two sides of this equality, it is reasonable to expect the following
\begin{conjecture}\label{conj 0}
If $f$ is relatively minimal and the singular fiber $f^{-1}(s_i)$ is semi-stable, then
$$12 a_i=\chi_2^i-3I_2^i.$$
\end{conjecture}
In this direction, Eriksson-Freixas-Mourougane~\cite{EFM16} has obtained interesting results for the asymptotic expansion of the BCOV metric. Recently, they informed us that they have solved this conjecture positively.

  On the other hand, by using the curvature properties of the moduli space of polarized Calabi-Yau manifolds and Yau's Schwarz lemma, we prove the following Arakelov type inequalities:
  \begin{theorem}\label{Arakelov inequality 0}
  1. for $n=3$, we have
\begin{itemize}
\item[(i)] $f$ is isotrivial $ \Leftrightarrow \deg\mathcal{H}_e^{3,0}=0\Leftrightarrow 12\sum_{i}a_i-\chi\deg\mathcal{H}_e^{3,0}=0$;
\item[(ii)] if $f$ is non-isotrivial, then
\begin{equation}\label{Arakelov inequality CY3 0}
 0< 2\deg\mathcal{H}_e^{3,0} \leq \sum_{i}a_i-\frac{\chi}{12}\deg\mathcal{H}_e^{3,0}\leq 2\pi(2g-2+s)[(\sqrt{h^{2,1}}+1)^2+1].
 \end{equation}
\end{itemize}
2. for $n=4$, we have
\begin{itemize}
\item[(i)] $f$ is isotrivial $ \Leftrightarrow \deg\mathcal{H}_e^{4,0}=0 \Leftrightarrow 2\deg\mathcal{H}_e^{4,0}+\deg\mathcal{H}_e^{3,1}=0$;
\item[(ii)] if $f$ is non-isotrivial, then
\begin{equation}
 0<\deg\mathcal{H}_e^{4,0} \leq 2\deg\mathcal{H}_e^{4,0}+\deg\mathcal{H}_e^{3,1}\leq \pi(2g-2+s)(h^{3,1}+4).
 \end{equation}\label{Arakelov inequality CY4 0}
\end{itemize}
\end{theorem}
In Section VI of~\cite{GGK09}, more refined Arakelov inequalities (not necessarily for Calabi-Yau manifolds) for $n=1,2,3$ is proved and their proof is Hodge-theoretic. It should be pointed out that our Schwarz lemma approach is also suited to proving general Arakelov inequality for (abstract) variation of Hodge structures by using the curvature computations in~\cite{GS69}.\par
  For $n=3$, we have the following bound on the Euler number $\chi$ of the general fiber of $f$:
\begin{corollary}\label{coro CY3 0}If $f \colon{\mathcal{X}}\to S$ is not isotrivial, then
\begin{equation}
\frac{\sum_{i}a_i -2\pi(2g-2+s)[(\sqrt{h^{2,1}}+1)^2+1]}{\deg\mathcal{H}_e^{3,0}}\leq \frac{\chi}{12}\leq \frac{\sum_{i}a_i}{\deg\mathcal{H}_e^{3,0}}-2.
 \end{equation}
\end{corollary}

 \begin{theorem}\label{Euler number inequality nond 0}If $f \colon{\mathcal{X}}\to S$ has no singular fibers and is not isotrivial.
\begin{itemize}
\item[(i)] for n=3, we have
\begin{equation}
\frac{-\pi(2g-2)[(\sqrt{h^{2,1}}+1)^2+1]}{\deg\mathcal{H}^{3,0}}\leq \frac{\chi}{24}\leq -1;
 \end{equation}
\item[(ii)] for n=4, we have
\begin{equation}\label{chi ineq CY4 0}
 1-\frac{\deg\mathcal{H}^{2,1}}{2\deg\mathcal{H}^{4,0}}\leq \frac{\chi}{24}\leq \frac{2\pi(2g-2)(h^{3,1}+4)-\deg\mathcal{H}^{2,1}}{2\deg\mathcal{H}^{4,0}}.
 \end{equation}
\end{itemize}
\end{theorem}
In the $n=3$ case, the Euler number is given by $\chi = 2(h^{1,1} - h^{2,1})$, so we see that a non-isotrivial family of Calabi-Yau threefolds with $h^{2,1}=1$ (mirror quintic type) must have singular fibers, as was shown in~\cite{GGK09}.
In the $n=4$ case, assume $\omega_{H^3}=c_1(\mathcal{H}^{2,1})=0$, then it was shown in~\cite{FL05} that if $\chi > 24(h^{3,1}+2)$ or $\chi < 24$, there exists no complete curve in $\mathcal{M}$, the moduli space of polarized Calabi-Yau manifolds. In particular, if $f$ has no singular fibers, then
$ 1 \leq\frac{\chi}{24}\leq h^{3,1}+2$. If, however, we assume the total space $\mathcal{X}$ has trivial canonical bundle, then it was shown in~\cite{TZ13} that $f$ must be isotrivial.\par

  We would like to remark that the results in this note are consequences of known results. Theorem \ref{main theorem 1 0} can be viewed as an integrated version of a theorem of Fang, Lu, Yoshikawa. Conjecture \ref{conj 0} is obtained by comparing Theorem \ref{main theorem 1 0} and a result of
Green, Griffiths, Kerr. The characterization of isotriviality in Theorem \ref{Arakelov inequality 0} is an interpretation of the vanishing of the integrals of the Weil-Petersson form and the Hodge form. The Arakelov inequality in Theorem \ref{Arakelov inequality 0} is a consequence of the known curvature bound of the Hodge metric and Yau's Schwartz Lemma. The bound on the Euler number
  is a reformulation of the Arakelov inequality. What we do in this paper is to bring together the works of two different groups, Fang-Lu-Yoshikawa on the one hand and Green-Griffiths-Kerr on the other. While Fang-Lu-Yoshikawa concentrate in local problems, Green-Griffiths-Kerr is more
interested in global results. By transforming from local to global, and from global to local, we build a connection between analytical
method in Fang-Lu-Yoshikawa and algebraic method in Green-Griffiths-Kerr. We find that their works can in fact help
each other and enable us to know more about the nature of BCOV invariants and degenerations of Calabi-Yau manifolds. We now know that the asymptotic values of BCOV invariants are essentially integers, in fact they reflect the topological and intersection data of singular fibers which is not clear in previous works. Also the Euler number bound provides partial evidence for a conjecture of Yau which states that there are at most finitely many topological types of Calabi-Yau manifolds in each dimension, compare~\cite{H90}.

  This paper is organized as follows. In Section 1, we introduce the notion of BCOV invariant; In Section 2, we recall the definitions and basic properties of the Weil-Petersson metric, generalized Hodge metrics and Hodge bundles; In Section 3, we show that the total singularities can be expressed by the sum of asymptotic values of BCOV invariants. In Section 4, by using Yau's Schwarz lemma, we prove Arakelov type inequalities and Euler number bound for Calabi-Yau family over a compact Riemann surface.

\vskip 1\baselineskip \textbf{Acknowledgements.}  The second author would
 like to thank Professors Carlos Simpson and Matt Kerr for useful discussions. We would like to thank the referees for carefully reading our manuscript and their valuable comments. We would also like to thank Eriksson, Freixas and Mourougane for their useful comments. The work of K. Liu is partially supported by the National Natural Science Foundation of China, Grant No. 11531012.

\vskip 1\baselineskip
      \end{section}
%%%%%%%%%%%%%%%%%%%%%%%%%%%%%%%%%%%%%%%%%%%%%%%%%%%%%%%%%%%%%%%%%%%%%%%%%%%%%%%%%%%%%%%%%%%%%%%%%%%%%%%%%%%%%%%%%%%%%%%%%%%%%%%%%%%%%%%%%%%%%%%%%%%%%%%%%%%%%%%%%%
%                                                                         Section 1
%%%%%%%%%%%%%%%%%%%%%%%%%%%%%%%%%%%%%%%%%%%%%%%%%%%%%%%%%%%%%%%%%%%%%%%%%%%%%%%%%%%%%%%%%%%%%%%%%%%%%%%%%%%%%%%%%%%%%%%%%%%%%%%%%%%%%%%%%%%%%%%%%%%%%%%%%%%%%%%%%%
\begin{section}{BCOV invariants}
    In this section, we introduce the notion of BCOV invariant. We will mainly follow the presentation in~\cite{Yo14}, see~\cite{FLY08} for more details.
%%%%%%%%%%%%%%%%%%%%%%%%%%%%%%%%%%%%%%%%%%%%%%%%%%%%%%%%%%%%%%%%%%%%%%%%%%%%%%%%%%%%%%%%%%%%%%%%%%%%%%%%%%%%%%%%%%%%%%%%%%%%%%%%%%%%%%%%%%%%%%%%%%%%%%%%%%%%%%%%%%
%                                                                         Section 1.1
%%%%%%%%%%%%%%%%%%%%%%%%%%%%%%%%%%%%%%%%%%%%%%%%%%%%%%%%%%%%%%%%%%%%%%%%%%%%%%%%%%%%%%%%%%%%%%%%%%%%%%%%%%%%%%%%%%%%%%%%%%%%%%%%%%%%%%%%%%%%%%%%%%%%%%%%%%%%%%%%%%
\begin{subsection}{Analytic torsion and BCOV torsion}
Let $(M,g)$ be a compact K\"ahler manifold of dimension $n$ with K\"ahler form $\omega$.
Let $\square_{p,q}=(\bar{\partial}+\bar{\partial}^{*})^{2}$ be the $\bar{\partial}$-Laplacian acting on $C^{\infty}$ $(p,q)$-forms on $M$ or equivalently $(0,q)$-forms on $M$ with values in $\Omega_{M}^{p}$, where $\Omega_{M}^{1}$ is the holomorphic cotangent bundle of $M$ and
$\Omega_{M}^{p}:=\Lambda^{p}\Omega_{M}^{1}$. Let \{$\lambda_j$\} be the eigenvalues of $\square_{p,q}$, then it is well-known
that $0=\lambda_0\leq\lambda_1\leq\lambda_2\leq\cdots\leq \lambda_j \longrightarrow +\infty$. The {\em spectral zeta function} of $\square_{p,q}$ is defined as
$$
\zeta_{p,q}(s):=\sum_{\lambda_j>0}\lambda_j^{-s},
$$
where the multiplicities of the eigenvalues are taken into account.
Then $\zeta_{p,q}(s)$ converges on the half-plane $\{s\in{ \mathbb{C}};\,\Re s>\dim M\}$, extends to a meromorphic function on ${ \mathbb{C}}$,
and is holomorphic at $s=0$. By Ray-Singer~\cite{RS73}, the {\em analytic torsion} of $(M,\Omega^{p}_{M})$ is the real number defined as
$$
\tau(M,\Omega_{M}^{p})
:=
\exp\{-\sum_{q\geq0}(-1)^{q}q\,\zeta'_{p,q}(0)\}.
$$
Note that $\tau(M,\Omega_{M}^{p})$ depends not only on the complex structure of $M$ but also on the metric $g$.
\par
In~\cite{BCOV93b}, Bershadsky-Cecotti-Ooguri-Vafa introduced the following combination of analytic torsions.

\begin{definition}
\label{def:BCOV:torsion}
The {\em BCOV torsion} of $(M,g)$ is the real number defined as
$$
T_{\rm BCOV}(M,g)
:=
\prod_{q\geq0}\tau(M,\Omega_{M}^{p})^{(-1)^{p}}
=
\exp\{-\sum_{p,q\geq0}(-1)^{p+q}pq\,\zeta'_{p,q}(0)\}.
$$
\end{definition}

If $\omega$ is the K\"ahler form of $g$, then we often write $T_{\rm BCOV}(M,\omega)$ for $T_{\rm BCOV}(M,g)$.
In general, $T_{\rm BCOV}(M,g)$ does depend on the choice of K\"ahler metric $g$ and hence is not a holomorphic invariant of $M$.
When $M$ is a Calabi-Yau threefold, it is possible to construct a holomorphic invariant of $M$ from $T_{\rm BCOV}(M,g)$ by multiplying a correction factor.
%%%%%%%%%%%%%%%%%%%%%%%%%%%%%%%%%%%%%%%%%%%%%%%%%%%%%%%%%%%%%%%%%%%%%%%%%%%%%%%%%%%%%%%%%%%%%%%%%%%%%%%%%%%%%%%%%%%%%%%%%%%%%%%%%%%%%%%%%%%%%%%%%%%%%%%%%%%%%%%%%%
%                                                                         Section 1.2
%%%%%%%%%%%%%%%%%%%%%%%%%%%%%%%%%%%%%%%%%%%%%%%%%%%%%%%%%%%%%%%%%%%%%%%%%%%%%%%%%%%%%%%%%%%%%%%%%%%%%%%%%%%%%%%%%%%%%%%%%%%%%%%%%%%%%%%%%%%%%%%%%%%%%%%%%%%%%%%%%%
\subsection
{Calabi-Yau threefolds and BCOV invariants}
\label{sect:1.2}

A compact connected K\"ahler manifold $X$ is {\em Calabi-Yau} if $H^q(X,{\mathcal O}_{X})=0$ for $0<q<\dim X$ and $K_{X}\cong{\mathcal O}_{X}$,
where $K_{X}$ is the canonical line bundle of $X$.
Let $X$ be a Calabi-Yau threefold. Let $g=\sum_{i,j}g_{i\bar{j}}\,dz_{i}\otimes d\bar{z}_{j}$ be a K\"ahler metric on $X$ and
let $\omega=\omega_{g}:=\sqrt{-1}\sum_{i,j}g_{i\bar{j}}dz_{i}\wedge d\bar{z}_{j}$ be the corresponding K\"ahler form. we define
$$
{\rm Vol}(X,\omega):=\frac{1}{(2\pi)^{3}}\int_{X}\frac{\omega^{3}}{3!}.
$$
The covolume of $H^{2}(X, \mathbb{Z})_{\rm free}:=H^{2}(X,{\mathbb{Z}})/{\rm Torsion}$ with respect to $[\omega]$ is defined as
$$
{\rm Vol}_{L^{2}}(H^{2}(X,{\mathbb{Z}}),[\omega])
:=
\det(\langle{ e}_{i},{ e}_{j}\rangle_{L^{2},[\omega]})_{1\leq i,j\leq b_{2}(X)}.
$$
Here $\{{ e}_{1},\ldots,{ e}_{b_{2}(X)}\}$ is a basis of $H^{2}(X,{\mathbb{Z}})_{\rm free}={\rm Im}\{H^{2}(X,{\mathbb{Z}})\to H^{2}(X,{\mathbb{R}})\}$ over ${\mathbb{R}}$
and $\langle\cdot,\cdot\rangle_{L^{2},[\omega]}$ is the inner product on $H^{2}(X,{\mathbb{R}})$ induced by integration of harmonic forms.

\begin{definition}
\label{def:Bott:Chern:correction}
For a Calabi-Yau threefold $X$ equipped with a K\"ahler form, define
$$
A(X,\omega)
:=
\exp\left[-\frac{1}{12}\int_{X}
\log\left(
\sqrt{-1}\frac{\eta\wedge\bar{\eta}}{\omega^{3}/3!}
\frac{{\rm Vol}(X,\omega)}{\|\eta\|_{L^{2}}^{2}}
\right)\,c_{3}(X,\omega)\right],
$$
where $c_{3}(X,\omega)$ is the top Chern form of $(X,\omega)$, $\eta\in H^{0}(X,K_{X})$ is a nowhere vanishing section and $\|\eta\|_{L^{2}}$ is its $L^{2}$-norm, i.e.,
$$
\|\eta\|_{L^{2}}
:=
\frac{1}{(2\pi)^{3}}\int_{X}\sqrt{-1}\eta\wedge\overline{\eta}.
$$
\end{definition}

Note that $A(X,\omega)$ is independent of the choice of $\eta$ and $A(X,\omega)=1$ if $\omega$ is Ricci-flat.

\begin{definition}
\label{def:BCOV:invariant}
The {\em BCOV invariant} of $X$ is the real number defined as
$$
\tau_{\rm BCOV}(X)
:=
{\rm Vol}(X,\omega)^{-3+\frac{\chi(X)}{12}}{\rm Vol}_{L^{2}}(H^{2}(X,{\bf Z}),[\omega])^{-1}T_{\rm BCOV}(X,\omega)\,A(X,\omega).
$$
\end{definition}
It turns out $\tau_{\rm BCOV}(X)$ is independent of the metric chosen.

\begin{theorem}\cite[Th.\,4.16]{FLY08}
\label{thm:FLY}
For a Calabi-Yau threefold $X$, $\tau_{\rm BCOV}(X)$ is independent of the choice of a K\"ahler form on $X$.
\end{theorem}
So we see that if $\pi\colon{\mathcal{X}}\to S$ is a smooth family of Calabi-Yau threefolds, then $\tau_{\rm BCOV}(t):=\tau_{\rm BCOV}(X_t), \forall t \in S,$ defines a smooth function on S, where the smoothness is proved in~\cite[Cor.\,3.9]{BGS88c}.
 \end{subsection}

 \end{section}
%%%%%%%%%%%%%%%%%%%%%%%%%%%%%%%%%%%%%%%%%%%%%%%%%%%%%%%%%%%%%%%%%%%%%%%%%%%%%%%%%%%%%%%%%%%%%%%%%%%%%%%%%%%%%%%%%%%%%%%%%%%%%%%%%%%%%%%%%%%%%%%%%%%%%%%%%%%%%%%%%%
%                                                                         Section 2
%%%%%%%%%%%%%%%%%%%%%%%%%%%%%%%%%%%%%%%%%%%%%%%%%%%%%%%%%%%%%%%%%%%%%%%%%%%%%%%%%%%%%%%%%%%%%%%%%%%%%%%%%%%%%%%%%%%%%%%%%%%%%%%%%%%%%%%%%%%%%%%%%%%%%%%%%%%%%%%%%%
\begin{section}{Weil-Petersson metric, generalized Hodge metrics and Hodge bundles}
Let $(X,[\omega])$ be a polarized Calabi-Yau manifold of dimension $n\geq 3$, that is, $X$ is a Calabi-Yau manifold and $[\omega]=c_1(L)\in H^2(X,\mathbb{Z})$ is the first Chern class of a positive line bundle $L$ on $X$. Let $\mathcal{M}$ be the (coase) moduli space of the polarized Calabi-Yau manifold $(X,[\omega])$. By Viehweg~\cite{Vie95}, $\mathcal{M}$ is quasi-projective. Locally, $\mathcal{M}$ is identified as a finite discrete quotient of the local versal deformation space $\rm Def$ (Kuranishi space) of $X$. By the Bogomolov-Tian-Todorov theorem~\cite{Tia87,Tod89}, the base space (which is a priori only a complex analytic space) $\rm Def$ of the Kuranishi family $$\pi: (\mathfrak{X}, X)\to \rm (Def,0)$$ is smooth. Indeed, $\rm Def$ is an open subset of the linear space $H^1(X,\Theta_X)$, where $\Theta_X$ is the holomorphic tangent bundle of $X$, so we may assume $\rm Def$ is contractible. Since $H^0(X,\Theta_X)\cong H^0(X,\Omega_X^{n-1})=0$, the Kuranishi family $\pi$ is universal. Define

\begin{equation}
 H^1(X,\Theta_X)_\omega :=\{\theta\in H^1(X,\Theta_X) \mid \theta\lrcorner \omega=0 \in H^2(X,\mathcal{O}_X) \},
 \end{equation}

and $\rm Def_\omega :=\rm Def \cap H^1(X,\Theta_X)_\omega$, then $\rm Def_\omega$ consists of those local deformations of $X$ preserving the polarization $[\omega]$. In fact, since $H^2(X,\mathcal{O}_X)=0$, we have $\rm Def_\omega=\rm Def$. We thus have the associated period mapping $\rm Def\to D$, where $D$ is the period domain, i.e. the classifying space of polarized (weight $k$) Hodge structure of $(X,[\omega])$. By Griffiths~\cite[Cor.\,3.6]{Gri68}, $\rm Def\to D$ is a holomorphic mapping, and is an immersion (Local Torelli) for $k=n$. We remark that $\mathcal{M}$ is an orbifold and is locally covered by $\rm Def$, so when we work with metrics, curvatures of $\mathcal{M}$, we can treat these notions on $\rm Def$ instead.\par
    Now we recall some natural metrics on $\mathcal{M}$, see~\cite{FL05} for more details.
%%%%%%%%%%%%%%%%%%%%%%%%%%%%%%%%%%%%%%%%%%%%%%%%%%%%%%%%%%%%%%%%%%%%%%%%%%%%%%%%%%%%%%%%%%%%%%%%%%%%%%%%%%%%%%%%%%%%%%%%%%%%%%%%%%%%%%%%%%%%%%%%%%%%%%%%%%%%%%%%%%
%                                                                         Section 2.1
%%%%%%%%%%%%%%%%%%%%%%%%%%%%%%%%%%%%%%%%%%%%%%%%%%%%%%%%%%%%%%%%%%%%%%%%%%%%%%%%%%%%%%%%%%%%%%%%%%%%%%%%%%%%%%%%%%%%%%%%%%%%%%%%%%%%%%%%%%%%%%%%%%%%%%%%%%%%%%%%%%
\begin{subsection}{Weil-Petersson metric}
Let $t\in \rm Def$, the Kodaira-Spencer map is now an isomorphism
\begin{equation}
\rho:
T_{t}\rm Def \xrightarrow[]{\cong} H^{1}(X_{t},\Theta_t),
\end{equation}
where $\Theta_t$ is the holomorphic tangent
bundle of $X_t$.%\footnote{If $t$ is a singular point of $\mathcal M$, then we should replace $U$ by $\hat U$, the local uniformization of $U$.}%

Let $(t_1,\cdots, t_m)$ be a local holomorphic
coordinate system of $\mathcal M$, we define a Hermitian inner product on $T_t\mathcal
M$ by
\[
\left(\frac{\pa}{\pa t_i},\frac{\pa}{\pa \bar t_j}
\right)_{WP}
=\int_{X_t} A_{i\bar\beta}^\alpha\cdot
\overline{A_{j\bar\delta}^\gamma} g^{\delta\bar\beta}
g_{\alpha\bar\gamma} dV_{X_t},
\]
where $A_i=A_{i\bar\beta}^\alpha\frac{\partial}{\partial
t_{\alpha}}\otimes d\bar t^\beta$, $(i=1,\cdots,m)$ are the harmonic
representation of
$\rho(\frac{\partial}{\partial t_i})$. This inner product on
each $T_t U$ for $t\in\mathcal  M$ gives
a Hermitian metric on the moduli space $\mathcal M$, which is
called the Weil-Petersson metric. Equipped with the
Weil-Petersson metric, $\mathcal M$ is a K\"ahler
orbifold.%~\footnote{A K\"ahler orbifold metric is a K\"ahler metric on the smooth part of the orbifold and lifts to an invariant K\"ahler metric on each local uniformization. See~\cite{Ruan} for details of  orbifolds and vector bundles over orbifolds.}%

Let $\Omega$ be a (nonzero) holomorphic $(n,0)$-form
on
$X_t$. Define $\Omega\lrcorner\rho(\frac{\pa}{\pa
t_i})$ to be the contraction of $\Omega$ and
$\rho(\frac{\pa}{\pa t_i})$. The Weil-Petersson metric can
 be re-written as (cf.~\cite{Tia87}):
\begin{equation}
\left(\frac{\pa}{\pa t_i},\frac{\pa}{\pa \bar t_j}
\right)_{WP}
=-\frac{\int_{X_t}\Omega\lrcorner\rho(\frac{\pa}{\pa
t_i})\wedge\overline{\Omega\lrcorner\rho(\frac{\pa}{\pa
t_j})} }{\int_{X_t}\Omega\wedge\bar\Omega}.
\end{equation}
\end{subsection}
%%%%%%%%%%%%%%%%%%%%%%%%%%%%%%%%%%%%%%%%%%%%%%%%%%%%%%%%%%%%%%%%%%%%%%%%%%%%%%%%%%%%%%%%%%%%%%%%%%%%%%%%%%%%%%%%%%%%%%%%%%%%%%%%%%%%%%%%%%%%%%%%%%%%%%%%%%%%%%%%%%
%                                                                         Section 2.2
%%%%%%%%%%%%%%%%%%%%%%%%%%%%%%%%%%%%%%%%%%%%%%%%%%%%%%%%%%%%%%%%%%%%%%%%%%%%%%%%%%%%%%%%%%%%%%%%%%%%%%%%%%%%%%%%%%%%%%%%%%%%%%%%%%%%%%%%%%%%%%%%%%%%%%%%%%%%%%%%%%
\begin{subsection}{Generalized Hodge metrics and Hodge bundles}\label{Generalized Hodge metrics and Hodge bundles}
Recall that, for all $0\leq k \leq n$, and $p+q=k$ there are natural holomorphic vector bundles $PR^q\pi_{*}\Omega_{\frak X/\rm Def}^p$, called Hodge bundles, on $\rm Def$ (hence on $\mathcal M$), whose fiber is
$$(PR^q\pi_{*}\Omega_{\frak X/\rm Def}^p)_t=PH^q(X_t,\Omega^p_{X_t}),$$
where $PH^q(X_t,\Omega^p_{X_t})$ is the primitive cohomology of $X_t$. By abuse of notation ,we will always use the same symbol $\mathcal{H}^{p,q}$ to denote the Hodge bundle on $\rm Def$ and on $\mathcal{M}$.  \par
    By differentiating harmonic representatives, we have a holomorphic bundle map
\begin{equation}\frac{\pa}{\pa t_i}: \mathcal{F}^{p}\rightarrow \mathcal{H}/\mathcal{F}^{p},
\end{equation}
where $\mathcal{F}^{p}=\mathcal{H}^{p,k-p}\oplus\mathcal{H}^{p+1,k-p-1}\oplus\cdots\oplus\mathcal{H}^{k,0}$, $\mathcal{H}:=PR^k\pi_{*}(\mathbb{C}).$ %, and for any local section $\alpha$ of $\mathcal{H}^{p,q}$, $\frac{\pa}{\pa t_i}(\alpha):=\alpha\lrcorner \rho(\frac{\pa}{\pa t_i})$
 In this way, we get a natural holomorphic bundle map
\begin{equation}\label{df}T(\rm Def)\rightarrow\underset{1\leq p\leq k} {\oplus}{\rm Hom}\,(\mathcal{F}^{p}, \mathcal{H}/\mathcal{F}^{p}).
\end{equation}
Note that, this bundle map is just the differential of the period mapping $\rm Def\to D$. There are natural metrics on the Hodge bundles $\mathcal{F}^{p}$, hence on each ${\rm Hom}\,(\mathcal{F}^{p}, \mathcal{H}/\mathcal{F}^{p})$, induced by the Riemann-Hodge bilinear relations. Let $h_{PH^k}$ be the pull back of the metric on $\underset{1\leq p\leq k} {\oplus}{\rm Hom}\,(\mathcal{F}^{p}, \mathcal{H}/\mathcal{F}^{p})$ by (\ref{df}), then $h_{PH^k}$ is semi-positive. We use $\omega_{PH^k}$ to denote the K\"ahler form of $h_{PH^k}$ for all $k\leq n$. And let
\begin{equation}\omega_{H^k}:=\omega_{PH^k}+\omega_{PH^{k-2}}+\cdots.
\end{equation}
We call both $\omega_{H^k}$ and $\omega_{PH^k}$ to be the
generalized Hodge metrics. We note that when $k=n$, $\omega_{PH^n}$ is a positive (1,1) form by local Torelli. It's just the pull back of the usual Hodge metric
on $D$ and we will call $\omega_{H}:=\omega_{PH^n}$ the Hodge metric. There are several interesting relations between these metrics and Hodge bundles:
\begin{lemma}\label{lem 1} Let $m=\dim \mathcal M = h^{n-1,1}$, then we have
\begin{itemize}
\item[(1)]
for $n=3$, $\omega_H=(m+3)\omega_{WP}+{\rm Ric}(\omega_{WP})$;
\item[(2)]
for $n=4$,  $\omega_H=(2m+4)\omega_{WP}+2{\rm Ric}(\omega_{WP})$;
\item[(3)]
$\omega_{PH^k}=\sum_{0\leq p\leq k}p c_1(\mathcal{H}^{p,q})$;
\item[(4)]
$\omega_{H^k}=\sum_{0\leq p\leq k}pc_{1}(R^{k-p}\pi_{*}\Omega^{p}_{\frak X/ \rm Def})$;
\item[(5)]
$\omega_{WP}=c_1(\mathcal{H}^{n,0})$.
\end{itemize}
\end{lemma}
\begin{proof}(1) is in~\cite[Thm.\,1.2]{Lu01}, see also~\cite[Thm.\,3.2]{Wang03} for an alternative proof; (2) is in~\cite[Coro.\,6.4]{LS04}; (3),(4) is in~\cite[Prop.\,2.8]{FL05}; (5) is in~\cite[Thm.\,2]{Tia87}.
\end{proof}
 Let us consider the Hodge metric $ \omega_{H}$ on $\mathcal M$ which is the pull-back of the homogeneous metric on the period domain $D$ by the period map. For Calabi-Yau threefolds, it can also be defined by (1) of Lemma \ref{lem 1}, and is shown, up to a constant, to be the pull back of the usual Hodge metric on period domain $D$. Hence these are equivalent definitions. The holomorphic sectional curvature of $\mathcal M$ under Weil-Petersson metric is not negative in general. The Hodge metric has better curvature property, indeed, the holomorphic sectional curvature is bounded from above by $\frac{-1}{(\sqrt{m}+1)^2+1}$ for Calabi-Yau threefolds. For Calabi-Yau fourfolds, the upper bound is $\frac{-1}{m+4}$. See~\cite{Lu01,LS04}.
\end{subsection}
%%%%%%%%%%%%%%%%%%%%%%%%%%%%%%%%%%%%%%%%%%%%%%%%%%%%%%%%%%%%%%%%%%%%%%%%%%%%%%%%%%%%%%%%%%%%%%%%%%%%%%%%%%%%%%%%%%%%%%%%%%%%%%%%%%%%%%%%%%%%%%%%%%%%%%%%%%%%%%%%%%
%                                                                         Section 2.3
%%%%%%%%%%%%%%%%%%%%%%%%%%%%%%%%%%%%%%%%%%%%%%%%%%%%%%%%%%%%%%%%%%%%%%%%%%%%%%%%%%%%%%%%%%%%%%%%%%%%%%%%%%%%%%%%%%%%%%%%%%%%%%%%%%%%%%%%%%%%%%%%%%%%%%%%%%%%%%%%%%
 \begin{subsection}{Weil-Petersson form and Hodge form}
Let $f: \mathcal{X} \to S$ be a smooth polarized family of Calabi-Yau manifolds with $(X_0:=f^{-1}(0), [\omega_0])\cong (X, [\omega])$, where $0\in S$ and $[\omega_0]$ is the polarization of $X_0$ induced from that of $\mathcal{X}$. Since $\mathcal{M}$ is the moduli space of the polarized Calabi-Yau manifold $(X, [\omega])$, we have a natural commutative diagram :
\[\label{commu diagram}\xymatrix{
   S  \ar[r]^\phi \ar[d]^\psi    & D/\Gamma  \ar[d]  \\
   \mathcal{M}  \ar[r]^\varphi    & D/G_\mathbb{Z}~~, }\]
where $\Gamma$ is the monodromy group of the family $f$, and $G_\mathbb{Z}=Aut(H_\mathbb{Z},Q)$, see~\cite{Gri84} for more information. To make notations simple, we will just use $\omega_{WP}$ and $\omega_{H}$ to denote the pull back form $\psi^*\omega_{WP}$ and $\psi^*\omega_{H}$, and call them the Weil-Petersson form and Hodge form on $S$, respectively. Similarly, we will use $\mathcal{H}^{p,q}$ to denote the pull back bundle $\psi^*\mathcal{H}^{p,q}$, which is isomorphic to $PR^qf_{*}\Omega_{\mathcal{X}/ S}^p$, and call them Hodge bundles on $S$.
\end{subsection}
\end{section}
%%%%%%%%%%%%%%%%%%%%%%%%%%%%%%%%%%%%%%%%%%%%%%%%%%%%%%%%%%%%%%%%%%%%%%%%%%%%%%%%%%%%%%%%%%%%%%%%%%%%%%%%%%%%%%%%%%%%%%%%%%%%%%%%%%%%%%%%%%%%%%%%%%%%%%%%%%%%%%%%%%
%                                                                         Section 3
%%%%%%%%%%%%%%%%%%%%%%%%%%%%%%%%%%%%%%%%%%%%%%%%%%%%%%%%%%%%%%%%%%%%%%%%%%%%%%%%%%%%%%%%%%%%%%%%%%%%%%%%%%%%%%%%%%%%%%%%%%%%%%%%%%%%%%%%%%%%%%%%%%%%%%%%%%%%%%%%%%
 \begin{section}{Total singularities in terms of asymptotic values of BCOV invariants}
Let $\mathcal{X}$ be a smooth projective variety of dimension $n+1$, $S$ be a compact Riemann surface and $f \colon{\mathcal{X}}\to S$ be a surjective, flat holomorphic map with generic fiber Calabi-Yau $n$-fold. Let $f^0 \colon{\mathcal{X}}^0 \to S^0$ be the smooth part of $f$, that is, each fiber $f^{-1}(s)$ is smooth for $s\in S^0$ and singular for $s\notin S^0$, $\mathcal{X}^0=f^{-1}(S^0)$ and $f^0$ is the restriction of $f$ to $\mathcal{X}^0$. We will use these notations throughout this paper unless otherwise stated. \par

Now, assume $n=3$, let $0\in S$ and $X_0=f^{-1}(0)$ be a singular fiber, $(U,t)$ be a coordinate neighborhood of $S$ centered at $0$. The following result is essential for us:
\begin{theorem}\cite[Th.\,10.1]{FLY08}
\label{thm:FLY 10.1}
Set
$$
a_0
:=
\lim_{t\to0}
\frac{\log\tau_{\rm BCOV}|_{U}(t)}
{\log|t|^{2}}
\in{\mathbb R}.
$$
Then the following equation of currents on $U$ holds:
\begin{equation}\label{current equation}
dd^{c}\log\tau_{\rm BCOV}
=
-\frac{\chi}{12}\,
\Omega_{\rm WP}
-\Omega_{\rm H}
+ a_0\,\delta_{0},
\end{equation}
 where $\tau_{\rm BCOV}(t):=\tau_{\rm BCOV}(X_t)$, for all $t\in S^0$, and $\log\tau_{\rm BCOV}$ is a locally integrable function on $S$.
  The currents $\Omega_{\rm WP}$ and $\Omega_{\rm H}$ are the trivial extensions of Weil-Petersson form and Hodge form from $S^0$ to $S$. Here $\delta_{0}$ is the Dirac current of the point $0$ and $\chi$ is the topological Euler number of a general fiber of $f$ .\end{theorem}

\begin{remark}\label{rk 3.2} In~\cite[Prop.\,7.3]{FLY08}, it is proved that the Weil-Petersson form $\omega_{\rm WP}$ and Hodge form $\omega_{\rm H}$ are bounded by the Poincar\'e metric on the punctured unit disc $\Delta^*$, this implies that $\omega_{\rm WP}$ and $\omega_{\rm H}$ are locally integrable on $S$ and extend trivially to closed positive (1,1)-currents on $S$. Recall that the trivial extension is defined as follows: if $\tilde{T}$ is the trivial extension of the (1,1)-form $T$ from $\Delta^*$ to the unit disc $\Delta$,
\begin{equation}\label{trivial extension} \tilde{T}(\eta):= \int_{\Delta}~\eta T,~~~~~~~~~\eta\in C_0^{\infty}(\Delta).
\end{equation}
\end{remark}

Recall that an n-dimensional singularity is an ordinary double point (ODP for short) if it is isomorphic to the hypersurface singularity at $0\in \mathbb{C}^{n+1}$ defined by the equation $z_0^2+z_1^2+\cdots + z_n^2=0$. And a fiber of $f \colon{\mathcal{X}}\to S$ is semi-stable if it is a reduced normal crossing divisor of $\mathcal{X}$.\par

   If $X_0$ only has ODP's as singularities, then it was shown by Yoshikawa~\cite[Th.\,5.2]{Yo14}, that $a_0=\frac{\#{\rm Sing}\,X_{0}}{6}$. Indeed,
   \begin{equation}\label{ODP}
   \log\tau_{\rm BCOV}(X_{t})=\frac{\#{\rm Sing}\,X_{0}}{6}\,\log|t|^{2}+O\left(\log(-\log|t|)\right), \qquad  t\to 0.
   \end{equation}
The value of $a_0$ for semi-stable degenerations is also obtained in~\cite[Th.\,3.11]{Yo14}, but has a very complicated expression. Inspired by the work of Green-Griffiths-Kerr~\cite{GGK09}, we conjecture that it should be expressed as the topological and intersecting data of the singular fiber, see Conjecture \ref{conj}. \par
   Our first observation is that the Weil-Petersson form and Hodge form and their trivial current extensions are naturally related to Hodge bundles. To begin with, let $\omega_{WP}, \, \omega_{H}$ be the Weil-Petersson form and Hodge form on $S^0$, respectively. $\Omega_{\rm WP}$ and $\Omega_{\rm H}$ are their trivial extensions on $S$ defined by (\ref{trivial extension}). On $S^0$, from Lemma \ref{lem 1} we know that

\begin{equation}\label{tian formula}\omega_{WP}=c_1(\mathcal{H}^{3,0}),\end{equation}
   and
\begin{equation}\label{FL formula}\omega_{H}=\omega_{PH^3}=\sum_{p=0}^3 pc_1(\mathcal{H}^{p,3-p})=3c_1(\mathcal{H}^{3,0})+c_1(\mathcal{H}^{2,1}).\end{equation}
   On the other hand, since we have assumed that all the local monodromy transformations around points $s_i\in S\setminus S^0$ are unipotent,
   there is then the canonical Deligne's extension $\mathcal{H}_e^{p,q}$ of the Hodge bundle $\mathcal{H}^{p,q}$ from $S^0$ to $S$. As is already remarked in~\cite[p.475]{GGK09},  it follows from Schmid's work~\cite{Sch73} that

\begin{itemize}
\item[(i)]
the forms $c_1(\mathcal{H}^{p,q})$ are integrable and define closed, $(1,1)$ currents $\overline{c_1(\mathcal{H}^{p,q})}$ (trivial extensions) on the completion $S$;
\item[(ii)]
$\overline{c_1(\mathcal{H}^{p,q})}= c_1(\mathcal{H}_e^{p,q})$ in $H^2_{DR}(S)$.
\end{itemize}
For a precise statement with a proof of (i) and (ii), we refer the reader to~\cite[Coro.\,5.23]{CKS86}.\par
We make the following observations:
\begin{lemma} \label{lem: currents of Hodge bundle}
 As cohomology classes in $H^2_{DR}(S)$, we have
\begin{equation}\label{l1}\Omega_{\rm WP}= \overline{c_1(\mathcal{H}^{3,0})}= c_1(\mathcal{H}_e^{3,0}),\end{equation} and
\begin{equation}\label{l2} \Omega_{\rm H}= 3\overline{c_1(\mathcal{H}^{3,0})}+\overline{c_1(\mathcal{H}^{2,1})}= 3c_1(\mathcal{H}_e^{3,0})+c_1(\mathcal{H}_e^{2,1}).\end{equation}
\end{lemma}
\begin{proof} This follows from (\ref{tian formula}),(\ref{FL formula}) and (ii). \end{proof}

We now have the following general picture:

\begin{theorem}[=Theorem \ref{main theorem 1 0}]\label{main theorem 1}
 Let $E:=S\setminus S^0=\{s_i\}\subset S$ be the discriminant locus of $f$ and $\tau_{\rm BCOV}(t)$ be the BCOV invariant of the fiber $f^{-1}(t)$. Set
 $$a_i
:=
\lim_{t\to s_i}
\frac{\log\tau_{\rm BCOV}(t)}
{\log|t|^{2}}.$$
Then
\begin{equation}\label{general picture}
\text{\rm Total Singularities}= (\chi+36)\deg\mathcal{H}_e^{3,0}+12 \deg\mathcal{H}_e^{2,1}=12\sum_{i}~a_i .
\end{equation}
\end{theorem}

\begin{proof} First note that by Theorem \ref{thm:FLY 10.1}, we have on $S$ the following current equation
$$dd^{c}\log\tau_{\rm BCOV}
=
-\frac{\chi(X)}{12}\,
\Omega_{\rm WP}
-\Omega_{\rm H}
+ \sum_{i}~a_i\,\delta_{s_i}.$$
By (\ref{l1}), (\ref{l2}) , we then have the following equality of cohomology classes in $H^2_{DR}(S)$:
$$dd^{c}\log\tau_{\rm BCOV}
=
-\frac{\chi +36}{12}\,
 c_1(\mathcal{H}_e^{3,0})
-c_1(\mathcal{H}_e^{2,1})
+ \sum_{i}~a_i\,\delta_{s_i}.$$

By integration on $S$, and noting that  $dd^{c}\log\tau_{\rm BCOV}$ is an exact current on $S$, we get
$$(\chi+36)\deg\mathcal{H}_e^{3,0}+12 \deg\mathcal{H}_e^{2,1}=12\sum_{i}~a_i.$$
\end{proof}

It was shown in \cite[Th.\,3.11]{Yo14}, that for each $i$, $a_i\in \mathbb{Q}$. Now from (\ref{general picture}), we see that
\begin{corollary}$$12\sum_{i}a_i\in \mathbb{Z} .$$
\end{corollary}

\begin{remark} By the monodromy theorem (see, e.g.,~\cite[p.41]{Gri84}), the local monodromy transformations $T_i$ around points $s_i$ in the discriminant locus $E=S\setminus S^0$ are quasi-unipotent, and already unipotent for semi-stable degenerations (see~\cite[p.185]{KK98},~\cite{Lan73}). In the general quasi-unipotent case, we may pass to the unipotent case by a sequence of base changes, or even better, to the semi-stable case by a sequence of semi-stable reductions, then if $\tilde{a_i}$ is the the asymptotic value for the semi-stable family, the asymptotic value for the original family is $a_i=\frac{\tilde{a_i}}{d_i}$, where $d_i$ is the degree of the base change corresponding to a neighborhood of $s_i$. See~\cite[Th.\,3.12]{Yo14}. And (\ref{general picture}) becomes
\begin{equation}
(\chi+36)\deg\mathcal{H}_e^{3,0}+12 \deg\mathcal{H}_e^{2,1}=12\sum_{i}~\tilde{a_i}=12\sum_{i}~a_id_i,
\end{equation}
where $\mathcal{H}_e^{3,0}$ and $\mathcal{H}_e^{2,1}$ are bundles on the base curve of the semi-stable family.
\end{remark}

    Now, assume $f$ has no singular fibers other than those with ODP's, denote these fibers by $X_j=f^{-1}(s_j), s_j\in E$ . And let $\Delta_j= \#{\rm Sing}\,X_{j}$ be the number of ODP's on $X_j$. By Theorem \ref{main theorem 1} and (\ref{ODP}), we immediately get

 \begin{corollary}If the singular fibers of $f$ has at most ODP singularities, then
 \begin{equation}\label{ODP formula}
 (\chi+36)\deg\mathcal{H}_e^{3,0}+12 \deg\mathcal{H}_e^{2,1}=12\sum_{j}~a_j=2\sum_{j}\Delta_j.
 \end{equation}
 \end{corollary}

   Furthermore, if we allow $f$ to have semi-stable fibers (i.e. reduced normal crossing divisors), and assume that $f$ is relatively minimal, then one of the main results in~\cite[Th.\,(V.A.1)(b)]{GGK09} says the following:
\begin{equation}\label{GGK main}
 (\chi+36)\deg\mathcal{H}_e^{3,0}+12 \deg\mathcal{H}_e^{2,1}=2~\sum_{j}\Delta_j+\sum_i\{\chi_2^i-3I_2^i\},
\end{equation}
where the second sum is over all semi-stable fibers $\{X_{s_i} \}$, and if $X_{s_i}=\bigcup_{\alpha} X_{\alpha}^i, X_{\alpha\beta}^i=X_{\alpha}^i\cap X_{\beta}^i$, then
$$\chi_2^i := \sum_{\alpha < \beta} \chi (X_{\alpha\beta}^i),$$
$$I_2^i := \sum_{\alpha < \beta} \deg [(X_{\alpha\beta}^i)^2]-\deg [(\sum_{\alpha < \beta}X_{\alpha\beta}^i)^2].$$
So we see that (\ref{ODP formula}) is a special case of (\ref{GGK main}) when there are only ODP singular fibers. This shouldn't be too surprising because our approach relies essentially on the curvature formula of determinant bundle~\cite[Th.\,0.1]{BGS88a}, which can be viewed as a differential form version of the \textit{Grothendieck-Riemann-Roch} formula (GRR for short). \par
    As another example, we consider a family $f: X\to B$ of curves of genus $g$ over a compact Riemann surface $B$ of genus $q$, assume the fibers of $f$ have at most ODP singularities, then a direct application of the GRR formula will give us:
\begin{equation}12\deg f_* \omega_{X/B}=\omega_{X/B}^2+e(X)-4(g-1)(q-1),
\end{equation}
where $e(X)$ is the topological Euler number of the algebraic surface $X$. Since $e(X)=4(g-1)(q-1)+\delta$, where $\delta$ is the sum of number of ODP's on each singular fiber, we then have
\begin{equation}\label{curve formula}12\deg f_* \omega_{X/B}=\omega_{X/B}^2+\delta.
\end{equation}
This formula and its variant are well-known to experts and can be deduced from the classical Riemann-Roch and Leray spectral sequence. The point is that it can be shown using analytic torsions. This is another manifestation of the equivalence of these two approaches we described above. In fact, let $\xi$ be a holomorphic vector bundle over $X$, $\lambda(\xi)$ be the determinant line bundle with fiber over $b\in B$ being $[\det H^0(X_b,\xi)]^{-1}\otimes \det H^1(X_b,\xi)$, in other words, $\lambda(\xi)=[\det f_*\xi]^{-1}\otimes \det R^1f_*\xi$, then a theorem~\cite[Th.\,2.1]{BB90} of Bismut and Bost says that
\begin{equation}c_1(\lambda(\xi),\parallel \parallel_Q)= -\int_f[Td(\omega_{X/B}^{-1})ch(\xi)]^{(4)}-\frac{1}{12}(rk \xi)\sum_i \Delta_i\delta_{s_i},
\end{equation}
where $\int_f$ means integration along a general fiber, $[]^{(4)}$ denotes the degree $4$ component, $\Delta_i=\#{\rm Sing}\,X_{s_i}$ is the number of ODP's on $X_{s_i}$ and $\delta_{s_i}$ is the Dirac current of the point $s_i$ as usual. Now, let $\xi=\mathcal{O}_X$, then
\begin{equation}\label{ls}
c_1(\det R^1f_*\mathcal{O}_X)= -\frac{1}{12}c_1^2(\omega_{X/B})-\frac{1}{12}\sum_i \Delta_i\delta_{s_i}.
\end{equation}
Note that $R^1f_*\mathcal{O}_X=(f_*\omega_{X/B})^\vee$. The integration of (\ref{ls}) over the base curve $B$ will give us exactly (\ref{curve formula}).\par

Let us go back to the study of Calabi-Yau family, we have seen that the ODP case of (\ref{GGK main}) can be deduced from the current equation (\ref{current equation}), conversely, we may use (\ref{GGK main}) to say something about the asymptotic values of the BCOV invariants in the semi-stable case, firstly we have:\par

 If $f$ is relatively minimal and has no singular fibers other than ODP's and semi-stable ones, then
\begin{equation}\label{ODP+semi-stable formula}
 12(\sum_{j}a_j+\sum_{i}a_i) = 2\sum_{j}\Delta_j+\sum_i\{\chi_2^i-3I_2^i\},
\end{equation}
 where we have used different indices $i,j$ to indicate different type of singular fibers.\par

In particular, we have the following
\begin{corollary}\label{semi-stable cor}
If $f$ is relatively minimal and has at most semi-stable singular fibers, then
\begin{equation}\label{semi-stable formula}
 12\sum_{i}a_i=\sum_i\{\chi_2^i-3I_2^i\}.
\end{equation}
\end{corollary}

Note that in the equality (\ref{semi-stable formula}), for each fixed $i$, $a_i$ is determined by the local data of $f$ near the point $s_i$, and $\chi_2^i-3I_2^i$ is determined by the fiber $X_{s_i}=f^{-1}(s_i)$, so it is reasonable to expect the following
\begin{conjecture}[=Conjecture \ref{conj 0}]\label{conj}
If $f$ is relatively minimal and the singular fiber $f^{-1}(s_i)$ is semi-stable, then
$$12 a_i=\chi_2^i-3I_2^i.$$
\end{conjecture}

Conversely, if we can prove conjecture \ref{conj}, then we can recover (\ref{GGK main}).

\end{section}
%%%%%%%%%%%%%%%%%%%%%%%%%%%%%%%%%%%%%%%%%%%%%%%%%%%%%%%%%%%%%%%%%%%%%%%%%%%%%%%%%%%%%%%%%%%%%%%%%%%%%%%%%%%%%%%%%%%%%%%%%%%%%%%%%%%%%%%%%%%%%%%%%%%%%%%%%%%%%%%%%%
%                                                                         Section 4
%%%%%%%%%%%%%%%%%%%%%%%%%%%%%%%%%%%%%%%%%%%%%%%%%%%%%%%%%%%%%%%%%%%%%%%%%%%%%%%%%%%%%%%%%%%%%%%%%%%%%%%%%%%%%%%%%%%%%%%%%%%%%%%%%%%%%%%%%%%%%%%%%%%%%%%%%%%%%%%%%%
\begin{section}{Arakelov inequality and Euler number bound for Calabi-Yau manifolds}
As in Section 3, let $\mathcal{X}$ be a smooth projective variety of dimension $n+1$, $S$ be a compact Riemann surface and $f \colon{\mathcal{X}}\to S$ be a surjective, flat holomorphic map with generic fiber Calabi-Yau $n$-fold. The smooth part $f^0:\mathcal{X}^0\to S^0$ of $f$ is a polarized family of Calabi-Yau $n$-folds with the polarization induced from that of $\mathcal{X}$. Remember that we always assume that all the local monodromy transformations around points $s_i\in S\setminus S^0$ are unipotent.
Following ideas in~\cite{Liu96b,Liu96a}, we prove Arakelov type inequalities by using curvature properties of the moduli space of Calabi-Yau manifolds and Yau's Schwarz lemma. As is explained in~\cite{Liu96a}, Yau's Schwarz lemma~\cite[Th.\,2']{Yau78a} is still applicable if the target manifold is an orbifold.\par
  There are lots of work on Arakelov inequality in the literature, see~\cite{GGK09,Pet00,VZ05} and the references therein.
%%%%%%%%%%%%%%%%%%%%%%%%%%%%%%%%%%%%%%%%%%%%%%%%%%%%%%%%%%%%%%%%%%%%%%%%%%%%%%%%%%%%%%%%%%%%%%%%%%%%%%%%%%%%%%%%%%%%%%%%%%%%%%%%%%%%%%%%%%%%%%%%%%%%%%%%%%%%%%%%%%
%                                                                         Section 4.1
%%%%%%%%%%%%%%%%%%%%%%%%%%%%%%%%%%%%%%%%%%%%%%%%%%%%%%%%%%%%%%%%%%%%%%%%%%%%%%%%%%%%%%%%%%%%%%%%%%%%%%%%%%%%%%%%%%%%%%%%%%%%%%%%%%%%%%%%%%%%%%%%%%%%%%%%%%%%%%%%%%
\begin{subsection}{General Calabi-Yau family}
    Let $E=S\setminus S^0$ be the discriminant locus, $s:=\#E$ be the number of singular fibers of $f$. We will use the fact: $\omega_H\geq 2\omega_{WP}$, which follows from the curvature property of Hodge bundles, see e.g.~\cite[Coro.\,2.10]{FL05}.

\begin{theorem}[=Theorem \ref{Arakelov inequality 0}]
1. for $n=3$, we have
\begin{itemize}
\item[(i)] $f$ is isotrivial $ \Leftrightarrow \deg\mathcal{H}_e^{3,0}=0\Leftrightarrow 12\sum_{i}a_i-\chi\deg\mathcal{H}_e^{3,0}=0$;
\item[(ii)] if $f$ is not isotrivial, then
\begin{equation}\label{Arakelov inequality CY3}
 0< 2\deg\mathcal{H}_e^{3,0} \leq \sum_{i}a_i-\frac{\chi}{12}\deg\mathcal{H}_e^{3,0}\leq 2\pi(2g-2+s)[(\sqrt{h^{2,1}}+1)^2+1].
 \end{equation}
\end{itemize}
2. for $n=4$, we have
\begin{itemize}
\item[(i)] $f$ is isotrivial $ \Leftrightarrow \deg\mathcal{H}_e^{4,0}=0\Leftrightarrow 2\deg\mathcal{H}_e^{4,0}+\deg\mathcal{H}_e^{3,1}=0$;
\item[(ii)] if $f$ is not isotrivial, then
\begin{equation}\label{Arakelov inequality CY4}
 0<\deg\mathcal{H}_e^{4,0} \leq 2\deg\mathcal{H}_e^{4,0}+\deg\mathcal{H}_e^{3,1}\leq \pi(2g-2+s)(h^{3,1}+4).
 \end{equation}
\end{itemize}
\end{theorem}
\begin{proof}1.(i) First, $f^0:\mathcal{X}^0\to S^0$ induces a holomorphic map from $S^0$ to the moduli space $\mathcal{M}$ of the polarized Calabi-Yau threefold $(X_0, [\omega])$: $$\psi: S^0\to \mathcal{M}, $$
where $0\in S^0$. Note that the Hodge form $\psi^*\omega_H$ on $S^0$ is semi-positive and locally integrable on $S$ (see Remark \ref{rk 3.2}), hence $\int_{S^0}\psi^*\omega_H= 0$ iff $\psi$ is constant which is equivalent to $f$ being isotrivial. Similarly,  $\int_{S^0}\psi^*\omega_{WP}= 0$ iff $f$ is isotrivial. But the Hodge form $\psi^*\omega_H$ and the Weil-Petersson form $\omega_{WP}$ have trivial current extensions on $S$ and by Lemma \ref{lem: currents of Hodge bundle} and Theorem \ref{main theorem 1}, we have

$$\int_{S^0}\psi^*\omega_H=\int_{S}\Omega_H=3\deg\mathcal{H}_e^{3,0}+\deg\mathcal{H}_e^{2,1}=\sum_{i}a_i-\frac{\chi}{12}\deg\mathcal{H}_e^{3,0},$$
$$\int_{S^0}\psi^*\omega_{WP}=\int_{S}\Omega_{WP}=\deg\mathcal{H}_e^{3,0},$$
so (i) follows.\par

(ii) By the uniformization theorem and Gauss-Bonnet theorem for punctured Riemann surfaces, $S^0$ admits a complete metric $\omega$ with constant curvature $K_1=+1, 0, -1$ when the Euler number $e(S^0)=2-2g-s$ is positive, zero or negative, respectively. Recall that the holomorphic sectional curvature of $\mathcal{M}$ with respect to the Hodge metric $\omega_H$ is bounded from above by $K_2=\frac{-1}{(\sqrt{h^{2,1}}+1)^2+1}$, so by Yau's Schwarz lemma, if $\psi$ is nonconstant, then $K_1=-1$ and
$$\psi^*\omega_H\leq \frac{K_1}{K_2}\omega=[(\sqrt{h^{2,1}}+1)^2+1]\omega, $$
by integration on $S^0$, we get
\begin{equation}\label{ineq 1}
\int_{S^0}\psi^*\omega_H\leq \frac{K_1}{K_2}\int_{S^0}\omega=2\pi(2g-2+s)[(\sqrt{h^{2,1}}+1)^2+1].
\end{equation}
Now substituting equalities in (i) into this inequality, and using
$$\int_{S^0}\psi^*\omega_H\geq 2\int_{S^0}\psi^*\omega_{WP}>0,$$ we get 1.(ii).\par
2. The proof is almost the same as $1.$ and is omitted. \end{proof}
\begin{remark}Note that for $n=3$, $$2\deg\mathcal{H}_e^{3,0}=\sum_{i}a_i-\frac{\chi}{12}\deg\mathcal{H}_e^{3,0}\Leftrightarrow \deg\mathcal{H}_e^{3,0}+\deg\mathcal{H}_e^{2,1}=0.$$
\end{remark}

 Now for $n=3$, we can give a bound on the Euler number $\chi$ of the general fiber of $f$:
\begin{corollary}\label{coro CY3}If $f \colon{\mathcal{X}}\to S$ is not isotrivial, then
\begin{equation}
\frac{\sum_{i}a_i -2\pi(2g-2+s)[(\sqrt{h^{2,1}}+1)^2+1]}{\deg\mathcal{H}_e^{3,0}}\leq \frac{\chi}{12}\leq \frac{\sum_{i}a_i}{\deg\mathcal{H}_e^{3,0}}-2.
 \end{equation}
\end{corollary}

\end{subsection}

\begin{subsection}{Non-degenerate Calabi-Yau family}
It is well-known that a non-isotrivial family of elliptic curves over a compact Riemann surface must have singular fibers. This phenomena does not hold for higher dimensional Calabi-Yau manifold, see~\cite{GGK09} for such examples. For non-degenerate Calabi-Yau family, we have the following bound for the Euler number $\chi$:
\begin{theorem}[=Theorem \ref{Euler number inequality nond 0}]If $f \colon{\mathcal{X}}\to S$ has no singular fibers and is not isotrivial.
\begin{itemize}
\item[(i)] for n=3, we have
\begin{equation}
\frac{-\pi(2g-2)[(\sqrt{h^{2,1}}+1)^2+1]}{\deg\mathcal{H}^{3,0}}\leq \frac{\chi}{24}\leq -1;
 \end{equation}
\item[(ii)] for n=4, we have
\begin{equation}\label{chi ineq CY4}
 1-\frac{\deg\mathcal{H}^{2,1}}{2\deg\mathcal{H}^{4,0}}\leq \frac{\chi}{24}\leq \frac{2\pi(2g-2)(h^{3,1}+4)-\deg\mathcal{H}^{2,1}}{2\deg\mathcal{H}^{4,0}}.
 \end{equation}
\end{itemize}
\end{theorem}

\begin{proof}
(i) is a direct consequence of Theorem \ref{coro CY3}. For (ii), first note that $f \colon{\mathcal{X}}\to S$ is now a polarized family of Calabi-Yau fourfolds with the polarization induced from that of $\mathcal{X}$. And recall Fang-Lu's formula~\cite[Th.\,1.1]{FL05}:
$$\sum_{i=1}^{i=n}(-1)^i\omega_{H^i}-\frac{\sqrt{-1}}{2}\partial\bar{\partial}\log T=\frac{\chi}{12}\omega_{WP},  $$
where $\omega_{H^i}$ is defined in Section \ref{Generalized Hodge metrics and Hodge bundles},
$T$ is the BCOV torsion with respect to the unique Ricci-flat metric~\cite{Yau78}. By Lemma \ref{lem 1},
\begin{equation}
\omega_{H^4}-\omega_{H^3}+\omega_{H^2}-\omega_{H^1}=4c_1(\mathcal{H}^{4,0)})+2c_1(\mathcal{H}^{3,1})-c_1(\mathcal{H}^{2,1}),
 \end{equation}
so by integration on $S$, we have
\begin{equation}
4\deg\mathcal{H}^{4,0}+2\deg\mathcal{H}^{3,1}-\deg\mathcal{H}^{2,1}=\frac{\chi}{12}\deg\mathcal{H}^{4,0}.
 \end{equation}
(\ref{chi ineq CY4}) follows from this and (\ref{Arakelov inequality CY4}).

\end{proof}

\end{subsection}

\end{section}

\end{document}